\documentclass[letterpaper, 10pt, conference]{ieeeconf}  % Comment this line out
\IEEEoverridecommandlockouts                              % This command is only
\usepackage{graphics} % for pdf, bitmapped graphics files
\usepackage{graphicx}
\usepackage{epsfig} % for postscript graphics files
\usepackage{times} % assumes new font selection scheme installed
\usepackage{balance}
\usepackage{color}
\usepackage{epstopdf}
\usepackage{amsmath}
\usepackage{amssymb}
\usepackage{bm}
\usepackage{booktabs}
\usepackage{mathrsfs}
\usepackage{theorem}
\usepackage{algorithm}
\usepackage{algorithmic}
\usepackage{tabularx}
\usepackage{indentfirst}

\newtheorem{theorem}{Theorem}

\newtheorem{lemma}{Lemma}

\newtheorem{assumption}{Assumption}

\newtheorem{definition}{Definition}

\newcommand{\Proof}{\noindent \textit{Proof.}$\;\;$}

\title{\LARGE \bf
Parallel Explicit Tube Model Predictive Control
}

\author{Kai Wang, Yuning Jiang, Juraj Oravec, Mario E. Villanueva, Boris Houska% <-this % stops a space
\thanks{K. Wang, Y. Jiang, M. E. Villanueva and B. Houska are with School of Information Science and Technology, ShanghaiTech University, Shanghai, China. 
{\tt\small\{wangkai1, jiangyn, meduardov,bhouska\}@shanghaitech.edu.cn}}
\thanks{J. Oravec is with Faculty of Chemical and Food Technology, Slovak University of Technology in Bratislava, Slovak Republic.
{\tt\small\{juraj.oravec@stuba.sk\}}}
}

\begin{document}

\maketitle
\thispagestyle{empty}
\pagestyle{empty}

%%%%%%%%%%%%%%%%%%%%%%%%%%%%%%%%%%%%%%%%%%%%%%%%%%%%%%%%%%%%%%%%%%%%%%%%%%%%%%%%
\begin{abstract}
This paper is about a parallel algorithm for tube-based model predictive control. The proposed control algorithm solves robust model predictive control problems suboptimally, while exploiting their structure. This is achieved by implementing a real-time algorithm that iterates between the evaluation of piecewise affine functions, corresponding to the parametric solution of small-scale robust MPC problems, and the online solution of structured equality constrained QPs. The performance of the associated real-time robust MPC controllers is illustrated by a numerical case study.  
\end{abstract}

%%%%%%%%%%%%%%%%%%%%%%%%%%%%%%%%%%%%%%%%%%%%%%%%%%%%%%%%%%%%%%%%%%%%%%%%%%%%%%%%
\section{Introduction}
During the past two decades, there have been many suggestions
on how to increase the robustness of nominal model predictive control
schemes by taking external disturbance models into account~\cite{Mayne2005}.
An in-depth review of the numerous approaches, for example, based on min-max
robust dynamic programming~\cite{Diehl2004}, scenario-tree
MPC~\cite{Scokaert1998,Engell2009}, semi-definite programming reformulations~\cite{Kothare1996}, uncertainty-affine feedback parameterizations~\cite{Goulart2006}, as well as modern Tube MPC formulations~\cite{Rakovic2012b,Villanueva2017} would certainly go beyond the scope of this paper. However, we refer to~\cite{Rakovic2012b}  and~\cite{Houska2019} for review articles of existing robust MPC approaches.

The focus of this paper is on the implementation of real-time algorithms for Tube MPC~\cite{Rakovic2012b}. Although one could argue that there already exist many efficient real-time algorithms for nominal MPC~\cite{Diehl2002,Houska2011,Mattingley2010,Zavala2009}, the same does not hold true for general Tube MPC approaches. Here, one is facing two numerical challenges: firstly, in Tube MPC, one needs to replace the predicted vector-valued state trajectory with a set-valued tube, and, secondly, the optimization variable of the robust MPC problem is in general a feedback law. In the past, several suggestions have been made on how to overcome these challenges. For example, rigid tube parameterizations~\cite{Mayne2005} use affine feedback laws together with polytopic tubes of constant cross-sections. Other tube parameterizations include so-called
homothetic~\cite{Rakovic2012,Rakovic2013} and elastic tubes~\cite{Rakovic2016}, which are typically based on polytopic sets, as well as ellipsoidal parameterizations~\cite{Villanueva2017}.

Many numerical solution methods for MPC use online optimization methods represent the feedback law implicitly, as the solution of a parametric optimization problem that may either be evaluated exactly or approximately~\cite{Diehl2002}. In contrast to this, explicit MPC methods attempt to move the computational burden of MPC into an offline routine that works out an explicit solution map of a parametric QP or LP~\cite{Bemporad2002}. Notice that both types of real-time MPC methods can---depending on the problem size---achieve run-times in the milli- and microsecond range~\cite{Mattingley2010,Houska2011}. The explicit MPC approach eventually outperforms the online approaches, as long
as the solution map consists of a not too large number of regions~\cite{Kvasnica2015}. However, the worst-case complexity of the explicit solution map grows exponentially with the number of constraints in the problem. Online solvers tend to perform better as soon as one attempts to solve MPC problems for larger systems.

In the context of Tube MPC, the question of whether to use explicit
or online MPC methods, needs to be addressed independently of the nominal
case. In fact, explicit robust MPC can be surprisingly efficient, as an empirical
observation is that robust controllers are often characterized by smaller feasible
sets---which can lead to a smaller number regions of the explicit solution
map~\cite{Bemporad2001,Kouramas2013}. Nevertheless, at the same time, it
must be clear that these robust Explicit MPC controllers are affected by the
same curse of  dimensionality as nominal Explicit MPC. In contrast to this,
online real-time Tube MPC algorithms~\cite{Zeilinger2009,Zeilinger2014}, have
the potential to scale-up to larger problems. For example, in~\cite{Hu2018}
a Tube MPC problem for a quadcopter with $10$ states has been implemented.

This paper introduces a real-time algorithms for solving Tube MPC problems for uncertain linear systems with polytopic constraints. The controller is based on a parallelizable model predictive control algorithm~\cite{Oravec2017,Jiang2018} and alternates between the evaluation of precomputed piecewise affine maps and solving an equality constrained quadratic program. The resulting algorithm yields a controller with recursive feasibility, constraint satisfaction and asymptotic stability guarantees in the presence of bounded disturbances.

Section~\ref{sec::problem} introduces the Tube MPC problem. Section~\ref{sec::algorithm} introduces a parallel Explicit MPC algorithm and establishes a robust asymptotic stability result. Section~\ref{sec::caseStudy} presents a comparisons between Explicit MPC, the proposed approach and conventional online solvers for a benchmark case-study. Section~\ref{sec::conclusion} concludes the paper.

\textbf{Notation} The sets of non-negative and positive integers are denoted by
$\mathbb N$ and $\mathbb N_{+}$, respectively. The Minkowski sum and Pontryagin difference of two sets $\mathbb Y,  \mathbb Z\subset\mathbb{R}^{n}$ is denoted by
\[
\begin{array}{rcl}
\mathbb Y\oplus \mathbb Z &=& \{ y + z \,|\, y\in\mathbb Y,\, z\in\mathbb Z \} \\
\text{and} \quad \mathbb Y\ominus\mathbb Z &=& \{ y \,|\, \{y\}\oplus\mathbb Z \subset\mathbb Y \}, 
\end{array}
\]
respectively. The sets of symmetric positive semidefinite and positive definite matrices are denoted by $\mathbb{S}^{n}_{+}$ and $\mathbb{S}^{n}_{++}$. 

%%%%%%%%%%%%%%%%%%%%%%%%%%%%%%%%%%%%%%%%%%%%%%%%%%%%%%%%%%%%%%%%%%%%%%%%%%%%%%%%
\section{Problem Statement}
\label{sec::problem}
% This section introduces basic definitions and the formulation of tube MPC problems.
% 
% \subsection{Uncertain discrete-time linear systems}

This paper considers uncertain control systems of the form
\begin{equation}
\label{eq::system}
x_{k+1} = Ax_k + Bu_k + w_k,
\end{equation}
with given matrices $A\in \mathbb{R}^{n_{\mathrm{x}}\times n_{\mathrm{x}}}$
and $B\in \mathbb{R}^{n_{\mathrm{x}}\times n_{\mathrm{u}}}$. Here, 
$x_k\in \mathbb{R}^{n_{\mathrm{x}}}$ and $u_k\in \mathbb{R}^{n_{\mathrm{u}}}$ 
denote the  state and control vectors at time $k$, while $w_k\in\mathbb{R}^{n_{\mathrm{x}}}$
is the disturbance vector. The disturbance is assumed
to take values in a compact set $\mathbb{W}\subset\mathbb{R}^{n_w}$. The states and
controls are required to satisfy hard constraints of the form
\begin{equation}
\label{constraint:origin}
\forall k\in\mathbb{N}, \quad x_k\in\mathbb{X}\quad\text{and}\quad u_k\in\mathbb{U}\;,
\end{equation}
for given closed set $\mathbb{X}\subset\mathbb{R}^{n_x}$ and compact set $\mathbb{U}\subset\mathbb{R}^{n_u}$.

\begin{assumption}	
\label{ass::sets}
The state constraint set $\mathbb{X}$ is a convex polyhedron, and the control constraint set $\mathbb{U}$ as well as the disturbance set $\mathbb{W}$ are convex polytopes. In addition, they all contain the origin in their interior.
\end{assumption}

\subsection{Rigid Robust Forward Invariant Tubes} 

Tube MPC controllers use the following definition~\cite{Mayne2005}.

\begin{definition}
\label{def::RFIT}
A sequence $X = ( X_0,\,  X_{1},\, \ldots )$ of compact sets is called a robust forward 
invariant tube for~\eqref{eq::system}, if there exists a feedback law $\mu:\mathbb{N}\times\mathbb{R}^{n_x}\to\mathbb{U}$ such that the state of the closed-loop system
\begin{equation*}
\forall k\in\mathbb{N},\quad x_{k+1} = Ax_{k} + B\mu(k,x_{k}) + w_{k}
\end{equation*}
satisfies $x_{k'}\in  X_{k'}$ whenever $x_{k}\in{X}_{k}$, for all $k'\geq k$, regardless of the disturbance sequence. 
\end{definition}

We write the state as $x_{k} = q_{k} + z_{k}$
and introduce the linear feedback law $ \mu(k,x_{k}) = v_{k} + K ( x_k - q_k) $.
Here, $q_{k}\in\mathbb{R}^{n_x}$ denotes the nominal (disturbance-free) component such that
\begin{equation}
\label{eq::nominal}
q_{k+1} = Aq_{k} + Bv_{k}
\end{equation}
with control input $v_k \in \mathbb{R}^{n_u}$. Thus, $z_{k}\in\mathbb{R}^{n_x}$
denotes the local error component satisfying
\begin{equation}
\label{eq::error}
z_{k+1} = (A + BK )z_{k} + w_{k} \;.
\end{equation}
Let $Z$ denote a pre-computed robust invariant set for the local error dynamics. A rigid RFITs is a set sequence of the form $X_{k} = \{ q_{k} \} \oplus Z$, together with associated control tubes
$ U_{k} = \{v_{k}\} \oplus KZ \subset \mathbb{R}^{n_u}$, such that $X$ is a robust forward invariant by construction.

\subsection{Tube-Based Model Predictive Control}

Rigid Tube MPC methods proceed by solving receding-horizon optimal 
control problems of the form\footnote{For polytopic sets $Y\subset \mathbb{R}^{n}$ and 
$V \subset \mathbb{R}^m$ together with a matrix $M \in\mathbb{R}^{n\times m}$, the relation
$\{y\}\oplus MV \subset Y \Longleftrightarrow y\in Y\ominus MV$ holds. }
\begin{equation}
\label{eq::tmpc}
\begin{alignedat}{2}
V(x_{0}) = &\min_{q,v} 
&& \sum^{N-1}_{k=0} \ell(q_{k},v_{k}) + m(q_{N}) \\
&\text{s.t.} && \left \{
\begin{alignedat}{1}
 \forall k &= 0,\ldots,N-1 \\
  q_{k+1} &=  A q_{k} + B v_k \\
  q_{k} &\in \mathbb{X} \ominus Z \; , \; v_{k} \in \mathbb{U} \ominus KZ \\
  q_{N} &\in \mathbb{X}_{T} \; , \; x_{0} \in \{q_0\} \oplus Z \;,
\end{alignedat}
\right. 
\end{alignedat}
\end{equation}
with $x_{0}\in\mathbb{R}^{n_x}$ denoting the current measurement.
The stage and terminal costs are given by
\begin{equation*}
\ell(q,v) = q^\intercal Q q + v^\intercal R v \quad\text{and}\quad 
m(q) = q^\intercal P q \;,
\end{equation*}
respectively, with $P,Q\in\mathbb{S}^{n_x}_{++}$ and $R\in\mathbb{S}^{n_u}_{++}$.  We denote the parametric solution map of~\eqref{eq::tmpc} by $q^\star(x_0)$ and $v^\star(x_0)$.

\begin{assumption}
\label{ass::cost}
The terminal constraint set $\mathbb{X}_{T}$ and the matrices $P\in\mathbb{R}^{n_x\times n_x}$, $K\in\mathbb{R}^{n_u\times n_x}$ are such that 
\begin{enumerate}
\item  the inclusions $(A+BK)\mathbb{X}_{T} \subset \mathbb X_{T}$, 
$\mathbb{X}_{T}  \subset \mathbb{X}\ominus Z$, and $K\mathbb{X}_T  \subset \mathbb{U} \ominus KZ$
hold.    
\item $m((A+BK)q) + \ell(q,Kq)  \leq m(q)$, for all $q\in \mathbb{X}_{T}$.
\end{enumerate}
\end{assumption}

Finally, the rigid tube MPC feedback law is given by
$$u^{\star}(x_0) = v^{\star}(x_0) + K ( x_0 - q_0^\star(x_0) ) \; . $$

\subsection{Recursive Feasibility and Asymptotic Stability}

As a consequence of Assumptions~\ref{ass::sets} and the positive definiteness
	of $P,Q$ and $R$,~\eqref{eq::tmpc} is a strictly convex quadratic program.
If Assumption~\ref{ass::cost} holds, $\mathbb{X}_T$ is a positive invariant set
for the nominal dynamics and $x_{0}\in \mathbb{X}_{T}$ implies
$(A+BK) q_{N}^\star(x_{0}) \in \mathbb{X}_{T}$.
Together with the invariance property of $Z$ this implies that the inclusion
\begin{equation*}
(A+BK) q_{N}^\star(x_{0}) + (A+BK)z + w \in \mathbb{X}_{T} \oplus Z
\end{equation*} 
holds for all $(z,w)\in Z \times \mathbb{W}$. Recursive feasibility of the 
Tube MPC scheme now follows from the inclusions $\mathbb{X}_{T}\oplus  Z \subset \mathbb{X}$ and 
$K\mathbb{X}_{T}\oplus K Z \subseteq \mathbb{U}$. Similarly,
the second part of Assumption~\ref{ass::cost} together with the 
invariance property of $ Z$ guarantee a strict decrease of the objective
along the closed-loop trajectory~\cite{Mayne2005}.

% In the next section we provide a parallel explicit algorithm for efficiently 
% solving~\eqref{eq::tmpc}. Before moving on, a remark must be made. The developments 
% below can be extended---with minor modifications---to other Tube MPC parameterizations
% such as homothetic~\cite{Rakovic2012}, (fully) parameterized~\cite{Rakovic2012} and 
% elastic Tube MPC formulations~\cite{Rakovic2016}. The 
% main requirement is for the underlying optimal control problem to be a parametric 
% strictly convex quadratic program.  

%%%%%%%%%%%%%%%%%%%%%%%%%%%%%%%%%%%%%%%%%%%%%%%%%%%%%%%%%%%%%%%%%%%%%%%%%%%%%%%%
\section{Real-Time Tube MPC}
\label{sec::algorithm}

In this section, we propose a real-time algorithm to approximately solve~\eqref{eq::tmpc}. 
First, let us introduce the vectors 
\begin{equation*}
y_k=[q_k^\intercal\,\,v_k^\intercal]^\intercal\;,\;k=0,...,N-1
\end{equation*}
and $y_N = q_N$ together with their associated constraint sets
\begin{equation*}
\begin{aligned}
\mathbb Y_{0} &= \left\{  y_{0} \ \middle| \
\begin{aligned}
&x_{0} \in \{ q_{0}\} \oplus  Z \\
&q_{0} \in \mathbb{X} \ominus  Z, \ v_{0} \in \mathbb{U} \ominus K Z  \\
&Aq_{0} + Bv_{0} \in \mathbb{X} \ominus  Z
\end{aligned}
\right\} \\[0.1cm]
\mathbb Y_{k} &= \left\{  y_{k} \ \middle| \
\begin{aligned}
&q_{k} \in \mathbb{X} \ominus  Z, \ v_{k} \in \mathbb{U} \ominus K Z  \\
&Aq_{k} + Bv_{k} \in \mathbb{X} \ominus  Z
\end{aligned}
\right\}\\
\end{aligned}
\end{equation*}
and $\mathbb{Y}_{N} = \{ y_{N} \ | \ q_{N}\in\mathbb{X}_{T} \}$, as well as the shorthand 
$$J(y) = \sum_{k=0}^{N} J_k(y_k)\;.$$ 
Now,~\eqref{eq::tmpc} is equivalent to
\begin{equation}
\label{eq::tmpc2}
\begin{alignedat}{2}
V(x_0) = &\min_{y} \ && J(y), \\
&\ \text{s.t.} &&\left\{
\begin{aligned}
& \forall k \in\{ 0,\dots,N-2,\} \\
& Dy_{k+1} - Cy_{k} = 0 \quad\mid \lambda_{k+1} \\
& y_{N} - Cy_{N-1} = 0  \quad\mid \lambda_{N} \\
& y_k\in\mathbb Y_k,y_{N-1}\in\mathbb Y_{N-1},y_N\in\mathbb Y_N.
\end{aligned}\right. 
\end{alignedat}
\end{equation}
with stage costs
\begin{equation*}
J_k(y_k) = \ell(q_{k},v_{k})\quad \text{for}\ k=0,1,...,N-1,
\end{equation*}
and $J_N(y_N) = m(q_{N})$ as well as matrices $C = [A\;B]$ and $D =[I\;0]$. 
Here, $\lambda$ denotes the Lagrangian multiplier of the dynamic equation.

\subsection{Parallel Tube-based MPC Algorithm}

\begin{algorithm}[tbp!]
	\caption{Real-time Parallel Robust MPC}
	\textbf{Initialization:} \\
	Initial guesses $y^1 = [y_0^1,\dots,y_N^1] $ and $\lambda^1 = [\lambda_1^1,\dots,\lambda_{N}^1]$. \\
	\textbf{Online:}
	\begin{enumerate}
		\item Wait for new measurement $x_0$ and compute 
	\begin{equation*}
	f^{1} = J(y^{1}) + J^\star(\lambda^1)\,.
	\end{equation*}
	Here, $J^{\star}(\lambda^1)$ denotes the convex conjugate of $J$,
	\begin{equation*}
	\begin{split}
	J^{\star}(\lambda) =\max_{y}& \ -J(y) + \sum^{N-1}_{k=1}(D^\intercal\lambda_{k}-C^\intercal\lambda_{k+1})y_{k}\\
	&\;\;-\lambda_1^\intercal Cy_{0} + \lambda_N^\intercal y_N\;.
	\end{split}
	\end{equation*}
	\textbf{If} $f^1\geq\gamma^2 x_{0}^\intercal Q x_0$, rescale
	\begin{equation*}
	y^1\leftarrow y^1\sqrt{\frac{\gamma^2x_{0}^\intercal Q x_0}{f^1} }\quad\text{and}\quad
	\lambda^1\leftarrow \lambda^1\sqrt{\frac{\gamma^2x_{0}^\intercal Q x_0}{f^1} }\;.
	\end{equation*}
	\item \textbf{For} $j = 1 \to \overline{m}$ \textbf{do}
	\end{enumerate}
	\begin{enumerate}
		\item[a)] Compute $\xi^{j} = ( \xi^{j}_{0},\xi^{j}_{1},\ldots,\xi^{j}_{N} )$ 
		using~\eqref{eq::decoupled}.
		\item[b)] Compute $(y^{j+1},\,\Delta^{j})$ using~\eqref{eq::coupling} and set
		\begin{equation*}
		\lambda^{j+1} \leftarrow \lambda^{j}+ \Delta^{j}\;.
		\end{equation*}
		\item[] 	\textbf{End}
	\end{enumerate}
	\begin{enumerate}
		\item[3)] Send the input $[0\;I]\xi_0^{\overline{m}}+K(x_0-D\xi_0^{\overline{m}})$ to the real process.
		\item[4)] Set $y^1 \leftarrow [y_1^{\overline{m}},\ldots,y_N^{\overline{m}},0]$, $\lambda^1\leftarrow [\lambda_2^{\overline{m}},\ldots,\lambda_{N}^{\overline{m}},0]$, go to Step 1.
	\end{enumerate}
	\label{alg:Parallel}
\end{algorithm}

Algorithm~1, computes an approximate solution of~\eqref{eq::tmpc2}.
The algorithm is based on the Augmented Lagrangian Alternating Direction
Inexact Newton method (ALADIN)~\cite{Houska2016}, tailored for solving 
MPC problems in real-time~\cite{Oravec2017,Jiang2018}. As discussed in~\cite{Jiang2018}, Step~1) rescales $y^1$ and $\lambda^1$ with parameter $\gamma>0$ satisfying the following assumption.
\begin{assumption}\label{ass::rescaling}
The constant $\gamma$ at Step~1) of Algorithm~1 satisfies\[
J(y^\star) + J^\star(\lambda^\star)\leq \gamma^2 x_0^\intercal Qx_0\;.\]
\end{assumption}

\vspace{0.2cm}
\noindent
Here, the optimal value $J(y^\star)$ and $J^\star(\lambda^\star)$ can be precomputed. This rescaling step prevents the shifted initialized guesses in Step 4) from being far away from the origin. Step~2) is the main step of Algorithm~1, which include two substeps. Step~2.a) solves an augmented Lagrangian optimization problem
\begin{equation}
\label{eq::ALprob}
\xi^j = \operatorname*{argmin}_{\xi\in\mathbb{Y}} \ J( \xi ) + (G^\intercal \lambda^j)^{\intercal}\xi
+ (\xi - y^j)^\intercal H (\xi - y^{j})\;.
\end{equation}
Here, $\mathbb{Y} = \mathbb{Y}_{0} \times \mathbb{Y}_{1} \times \ldots \times \mathbb{Y}_{N}$,  
\begin{equation*}
H = \nabla^2 J(y)\;,\;G = \begin{pmatrix}
-C &  D  &         &        & 0    \\
   & -C  & D       &        &      \\
   &     & \ddots  & \ddots &      \\
0  &     &         & -C     & I        
\end{pmatrix}\;.
\end{equation*}
Notice that~\eqref{eq::ALprob} has a completely separable 
structure and 
$\xi^{j} = ( \xi_{0}^{j}, \xi_{1}^{j}\ldots, \xi^{j}_{N} )$
can be computed via
\begin{align}\notag
\xi^{j}_0 &= \operatorname*{argmin}_{\xi \in\mathbb{Y}_{0}} \ 
J_{0}(\xi) - (C^\intercal \lambda^{j}_{1} )^\intercal \xi + J_{0}( \xi - y_{0}^{j}) \\\notag
\xi^{j}_k &= \operatorname*{argmin}_{\xi\in\mathbb{Y}_{k}} \ 
J_{k}(\xi) + (D^\intercal \lambda^{j}_{k} - C^\intercal\lambda^{j}_{k+1} )^\intercal \xi + J_{k}( \xi - y_{k}^{j}) \\
\xi^{j}_N &= \operatorname*{argmin}_{\xi \in\mathbb{Y}_{N}} \ 
J_{N}(\xi) + ( \lambda^{j}_{N} )^\intercal \xi + J_{N}( \xi - y_{N}^{j}) \;.\label{eq::decoupled}
\end{align}
with $k\in\{1,\ldots,N-1\}$. 

In Step 2.b), the next iterate for the primal variables
and the increment for the dual variables is obtained by solving 
\begin{equation}
\label{eq::coupling}
\begin{alignedat}{2}
y^{j+1} = &\operatorname*{argmin}_{y} \ &&\sum^{N}_{k=0} J_{k}(y_{k}-2\xi^{j}_{k}+y^{j}_{k}) \\
& \ \ \text{s.t.} && \left\{
\begin{aligned}
&\forall k\in\{0,\ldots,N-2\} \\
&Dy_{k+1} - Cy_{k} = 0 \quad |\ \Delta_{k+1}^{j} \\
&y_{N} - Cy_{N-1} = 0 \quad |\ \Delta_{N}^{j} \;.
\end{aligned}\right.
\end{alignedat}
\end{equation}
This coupled QP can be interpreted as an unconstrained Linear Quadratic Regulator (LQR) problem tracking a weighted average of the solution of~\eqref{eq::ALprob} 
and the previous iterate $y^{j}$~\cite{Jiang2018}. In the following, we denote the suboptimal solutions from Step~2) by $q_0^\circ(x_0)=[I\;0]\xi^{\overline{m}}_{0}$ and $v_0^\circ(x_0)=[0\;I]\xi^{\overline{m}}_{0}$.

\subsection{Recursive Feasibility}
Despite the fact that the rigid Tube MPC controller is recursively feasible by design, one may ask if this property could get lost if~\eqref{eq::tmpc} is not solved to optimality. As it turns out, a suboptimal solution computed with Algorithm~1 preserves this property as long as the following assumption holds.
\begin{assumption}	
	\label{ass::control_invariant}
	We assume that the state constraint set $\mathbb{X}$ is robust control invariant.
\end{assumption}
This is a result of the construction of the constraint sets $\mathbb{Y}_{k}$ used in the decoupled problems. In particular, as the constraint $Aq_0 + Bv_0 \in \mathbb{X}\ominus  Z$ is enforced and 
$$v_0^\circ(x_0)+K(x_{0}-q_0^\circ(x_0))$$ is sent to the process,  we have that the closed-loop system satisfies 
\begin{equation}\label{eq::RTI_system}
x^{+}_{0} = A x_{0} + Bv_0^\circ(x_0)+BK(x_{0}-q_0^\circ(x_0)) + w \in \mathbb{X}
\end{equation}
for all $w\in\mathbb{W}$. Using the invariance properties of $ Z$ and $\mathbb{X}_{T}$ and Assumption~\ref{ass::control_invariant},
a recursive feasibility argument can be constructed along the lines of~\cite{Mayne2005}.

\subsection{Asymptotic Stability of Algorithm~1}
This section analyzes robust stability of the proposed closed loop scheme. In Proposition~3 in~\cite{Mayne2005} it has been shown that 
the value function $V$ satisfies the inequality,
\begin{equation}
\label{eq::standard_stability}
V(x_1^\star) \leq V(x_0) -J_0(y_0^\star(x_0))
\end{equation}
with
\begin{equation}\notag
\label{eq::standard_update}
x_1^\star = Ax_0 + B(v_0^\star(x_0) + K(x_0 - q_0^\star(x_0))) + w_0
\end{equation}
as long as Assumption~\ref{ass::sets} and~\ref{ass::cost} hold and Problem~\eqref{eq::tmpc} is feasible for the initial state $x_0$. In order to simplify the notations, we use $y^\star$ to denote $y^\star(x_0)$.  Now, we have 
\begin{equation}
\label{eq::our_stability}
V(x_0^+) \leq V(x_0)-(J_0(y_0^\star)-V(x_0^+)+V(x_1^\star))
\end{equation}
with
\begin{equation}\notag
\label{eq::our_update}
x^+_0=Ax_0+B(v_0^\circ(x_0)+K(x_0-q_0^\circ(x_0))+w_0\,,
\end{equation}
that is, $V$ is a Lyapunov function as long as
\begin{equation}
\label{eq::stability_condition}
J_0(y_0^\star)-V(x_0^+)+V(x_1^\star) \geq \alpha J_0(y_0^\star)\,,
\end{equation}
for a constant $\alpha>0$. Next, we use that
\begin{equation}
\label{eq::Algorithm_contraction}
\begin{split}
&J(y^{j+1}-y^\star) + J^\star(\lambda^{j+1}-\lambda^\star)\\
\leq\; & \kappa \left(J(y^{j}-y^\star) + J^\star(\lambda^{j}-\lambda^\star)
\right)
\end{split} 
\end{equation}
with a constant  $0<\kappa<1$, which has been shown in Theorem~1 of~\cite{Jiang2018}. A proof of the following
lemma can be found in~\cite{Jiang2018}, too.
\begin{lemma}\label{lem::eta_and_tau}
Let Assumption~\ref{ass::sets} hold, Problem~\eqref{eq::tmpc} is a strongly convex parametric QP such that the value function $V$ satisfies
\begin{equation}
	\label{eq::eta_and_tau}
	|V(x_0^+) - V(x_1^\star)| \leq \eta \|x_0^+-x_1^\star\|_Q + \frac{\tau}{2}\|x_0^+-x_1^\star\|_Q^2\;.
\end{equation}
with $\eta,\tau>0$.

% and the optimal solution $y_0^\star$ satisfies
%\begin{equation}
%	\label{eq::beta}
%	\|x_0\|_Q^2\leq \beta^2 J_0(y_0^\star)\;.
%\end{equation}	
%with $\beta>0$.
\end{lemma}

Now, the main idea for deriving a stability statement for Algorithm~1 is to show that the real-time approximation
$$x^\star_1 \approx x_0^+$$
of the optimal next state is sufficiently accurate to still ensure descent of the Lyapunov function. In order to bound the corresponding error term, the following technical result is needed.

\begin{lemma}\label{lem::stability_condition}
Let Assumption~\ref{ass::sets},~\ref{ass::cost} and~\ref{ass::rescaling} hold. There exists a constant $\sigma > 0$ such that
the iterate $x_0^+$ satisfies the inequality
\begin{equation}
\label{eq::iterate_bound}
\|x_0^+-x^\star_1\|_Q^2 \leq \sigma \kappa^{\overline{m}} J_0(y_0^\star)\;.
\end{equation}
\end{lemma}

\Proof Firstly, the equation
\begin{equation}
\label{eq::iteration_error}
x^+_0-x^\star_1 = \mathcal P (\xi^{\overline{m}} - y^\star)
\end{equation}
holds with $\mathcal P = [ -BK, \, B, \, 0, \ldots, 0 ]$. Because the  Algorithm~1 converges
globally with linear rate~\cite{Jiang2018}, there exists a constant $\tilde \sigma > 0$
such that
\[
\Vert \xi^{\overline{m}} - y^\star \Vert^2 \leq \tilde \sigma \kappa^{\overline m} \left( J(y^1-y^\star) 
+J^\star(\lambda^1-\lambda^{\star}) \right) \; .
\]
Now, Assumption~\ref{ass::rescaling}, the rescaling step in Algorithm~1, and~\eqref{eq::iteration_error}
imply that there must exist a constant $\sigma > 0$ such that
\begin{equation}
\|x_0^+-x_1\|_Q^2 \leq \sigma \kappa^{\overline{m}} J_0(y_0^\star) \; ,
\end{equation}
which is the statement of this lemma.\hfill$\square$\\[0.2cm]

\noindent
The main stability properties of Algorithm~1 can now be summarized as follows.

\begin{theorem}\label{Theory::stable2}
	Let Assumption~\ref{ass::sets},~\ref{ass::cost} and~\ref{ass::rescaling} hold, if the number of inner loops in Step~2) of Algorithm~\ref{alg:Parallel} satisfies
	\begin{equation}\label{eq::m}
	\overline{m}\geq 2\frac{\log \left(\eta \sqrt{\sigma} + \frac{\tau \sigma}{2}\right)}{\log(1/k)}\;,
	\end{equation}
	Algorithm~1 yields an asymptotically stable closed-loop controller. 
\end{theorem}
\Proof From inequalities~\eqref{eq::eta_and_tau} and~\eqref{eq::iterate_bound} in Lemma~\ref{lem::stability_condition}, we have
\begin{equation}
\label{eq::iterate_bound2}
\begin{split}
&|V(x_0^+) - V(x_1^\star)|\\
\leq\;&\left[\eta \sqrt{\sigma} + \frac{\tau \sigma}{2}\right]\kappa^{\frac{\overline{m}}{2}}J(y^\star_0)\;.
\end{split}
\end{equation}
Now, the inequality~\eqref{eq::m} in Theorem~\ref{Theory::stable2} follows directly from~\eqref{eq::iterate_bound2}.
By combining Lyapunov descent condition~\eqref{eq::stability_condition} and~\eqref{eq::iterate_bound2}, we have that $V$ can be used as a Lyapunov function that proves local asymptotic stability~\cite{Rawlings2017} with 
\[
\alpha = 1-\left[\eta \sqrt{\sigma} + \frac{\tau \sigma}{2}\right]\kappa^{\frac{\overline{m}}{2}}>0\;.
\]
\hfill$\square$\\[0.2cm]
\noindent
Notice that the statement of the above theorem implies that the set $Z$ is robustly stable. This result follows immediately from Theorem~1 in~\cite{Mayne2005}.
	
%%%%%%%%%%%%%%%%%%%%%%%%%%%%%%%%%%%%%%%%%%%%%%%%%%%

\section{Numerical Case Study}
\label{sec::caseStudy}

%%%%%%%%%%%%%%%%%%%%%%%%%%%%%%%%%%%%%%%%%%%%%%%%%%%%%%
We consider an uncertain control system of the form~\cite{Mayne2005}
\begin{equation}\notag
x_{k+1} =
\begin{pmatrix}
1      & 1      \\
0      & 1
\end{pmatrix} x_{k} +
\begin{pmatrix}
0.5 \\
1     
\end{pmatrix} u_{k} + w_{k}\;,
\end{equation}
with initial value $x_{0}^\intercal = -(7,2)$.
Here, the disturbance is assumed to take values on the set 
$\mathbb{W} = \{w\ |\  \|w\|_\infty\leqslant 0.1\}$. The state and control 
constraints are given by the sets 
\begin{equation*}
\mathbb{X} = \left\{ x \ \middle|\  
[0\;1]x \leq 2 
\right\} \quad\text{and}\quad
\mathbb{U} = \{u\ |\  |u|\leq 1 \},\\
\end{equation*}
The matrices for the stage and terminal costs are given by
$Q = I$, $R = 0.1$, $K=-( 0.62,\, 1.27 )$ and
\begin{equation*}
P = 
\begin{pmatrix}
2.06      & 0.60      \\
0.60      & 1.40
\end{pmatrix}\;.
\end{equation*}
Here, $P$ and $K$ were computed as the solution of an algebraic Riccati 
equation yielding the optimal LQR controller for the nominal system. The set $Z$
was computed using Algorithm~1 in~\cite{Rakovic2005} such that
\begin{equation*}
 Z_{\infty} \subseteq Z \subseteq  Z_{\infty} \oplus \{ x \ | \ \Vert x \Vert_{\infty}\leq 10^{-4}  \} \; ,
\end{equation*}
where $ Z_\infty$ denotes the minimal robust positive invariant set for the local
error dynamics~\eqref{eq::error}.
Likewise, the set $\mathbb{X}_{T}$ was chosen as the maximal positively invariant set
for the nominal dynamics~\eqref{eq::nominal}, and computed according to Algorithm~3.2 in~\cite{Gilbert1991}. It is easy to check that with this construction, Assumption~2 is satisfied. 

In order to perform comparisons, rigid Tube-based MPC controllers with horizons $N=10,20,\ldots,100$ were implemented based on solving~\eqref{eq::tmpc} 1) explicitly, 2) online using a centralized QP solver, and 3) online using Algorithm~1.  All the algorithms in this section were implemented in {\tt MATLAB} R2019a on a Windows 10 personal computer with an i7 3.6 GHz and 16 GB of RAM. Polyhedral computations and the explicit solution of quadratic programs were done using the 
Multi-Parametric Toolbox ({\tt MPT} v3.2.1)~\cite{Herceg2013} through the {\tt YALMIP} interface 
(R20181012)~\cite{Lofberg2004}. The online centralized Tube MPC controller was implemented 
using the {\tt qpOASES} solver~\cite{Ferreau2014} with hot-start option after condensing.

The following table shows the number of regions for the solution of problem~\eqref{eq::tmpc} for increasing time horizons. 
\begin{table}[htbp!]
	\centering
	\begin{tabular}{ccc}
		\toprule  
		$N$ &number of regions &memory [KB] \\ [0.1cm]
		\midrule  
		%5  & 611    & 38 \\
		10 & 1648   & 174\\[0.1cm]
		20 & 5312   & 1028\\[0.1cm]
		30 & 11050   & 3108\\[0.1cm]
		50 & 25160   & 11500\\[0.1cm]
		70 & 42700  & 27066\\
		\bottomrule  
	\end{tabular}
	\label{table::regions}
\end{table}
As expected, both the number of regions as well as the memory 
requirements grow very fast as $N$ increases. On the other hand, the memory requirements needed 
for the solution of~\eqref{eq::decoupled} in Algorithm~1 is 36 kB (corresponding to 465
critical regions), irrespective of the horizon length $N$.  
\begin{figure}[htbp!]
	\centering	  
	\includegraphics[clip,width=0.9\linewidth,trim={3.3cm 8cm 3.3cm 8.3cm}]{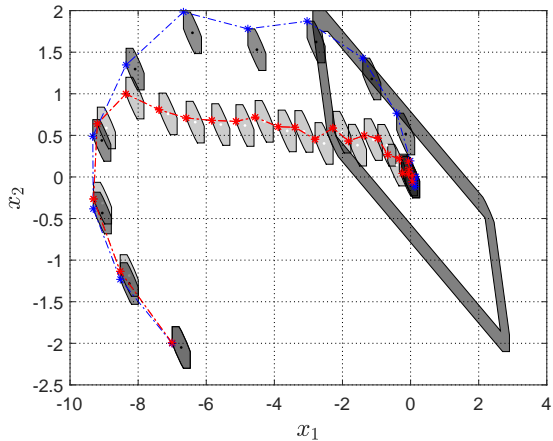}  
	\caption{Closed-loop simulation for rigid Tube MPC ($N=20$) scheme using Algorithm~1. 		Rigid RFITs computed with $\overline{m} = 2$ are shown in light gray, while those computed with $\overline{m}=5$ are shown in darker gray. The closed loop trajectories computed with $\overline{m}=2$ and 
		$\overline{m}=5$ are shown in red and blue, respectively. The terminal set $\mathbb{X}_{T}$ is depicted in white and $\mathbb{X}_{T}\oplus  Z$ in dark gray.}
	\label{fig::tube} 
\end{figure}

Figure~\ref{fig::tube} shows a closed-loop simulation based on the rigid tube MPC with Algorithm~1 (with $N=20$) for $\overline{m} = 2$ (red line for state, light gray for the tube) and $\overline{m} = 5$ (blue line for state, darker gray for the tube).
% Notice that both controllers drive the system to $Z$ in an asymptotically stable manner. 

\begin{figure}[htbp!]
	\centering  
	\includegraphics[clip,width=0.9\linewidth,trim={0.3cm 0.27cm 0.3cm 0.45cm}]{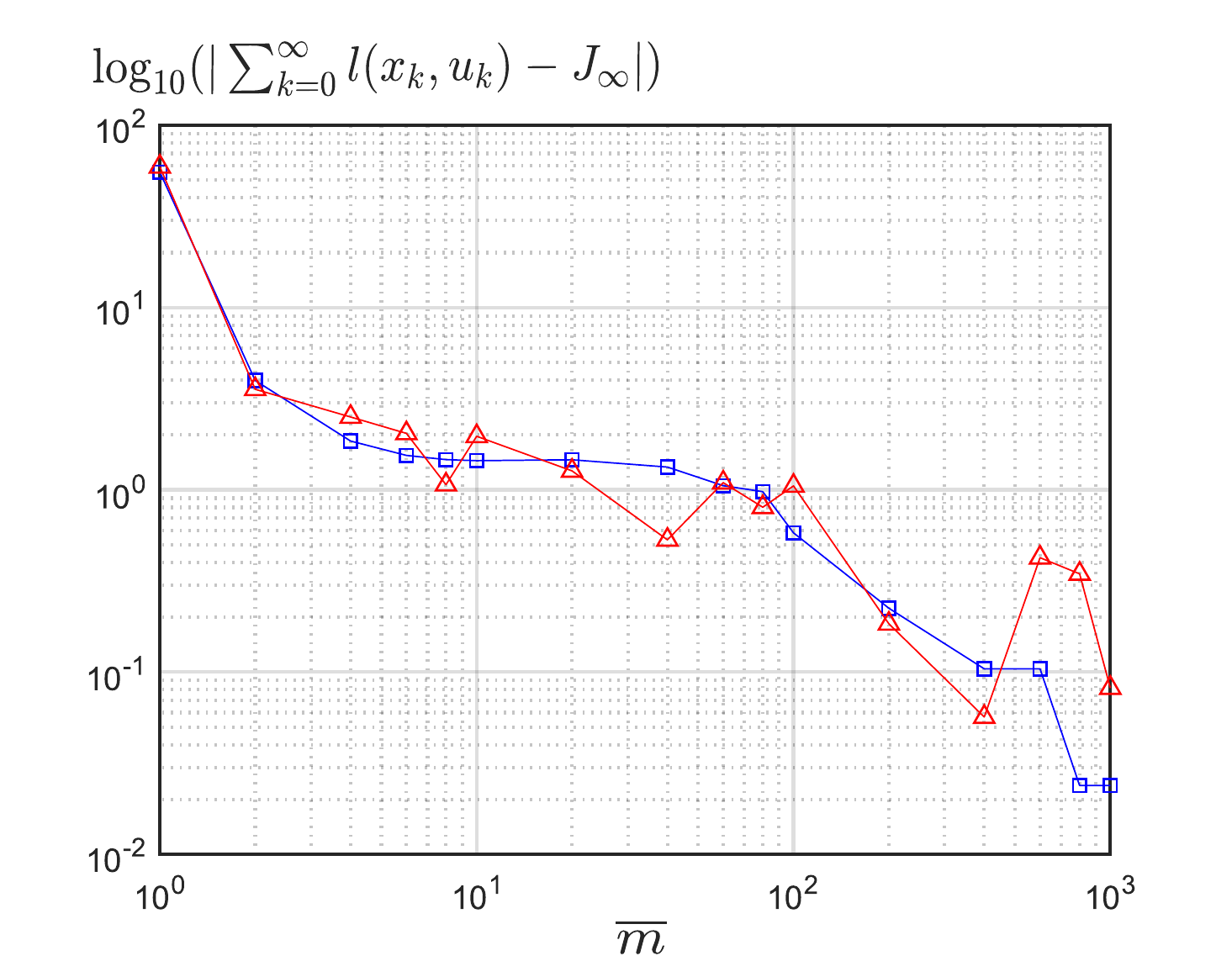}  
%		\put(0,20){\rotatebox{90}{$\log{\left| \sum^{\infty}_{k=0} \ell(x_{k},u_{k}) - J_{\infty} \right|}$}}
%		\put(55,0){$\overline{m}$}
	\caption{Closed-loop performance degradation for Algorithm~1 with respect to the optimal cost $J_{\infty}$.
		The blue line with square markers denotes the performance loss for $w = 0$, while the red line with triangular
		markers denotes the same for process noise.}  
	\label{fig::performance} 
\end{figure}
 
Figure~\ref{fig::performance} shows the performance loss for the same controller with
respect to the optimal cost $J_\infty$ as $\overline{m}$ increases. It is easy to see that 
in the absence of uncertainty (blue line with square markers), the performance loss tends
to zero.
\begin{figure}
	\centering  
	\includegraphics[width=0.9\linewidth]{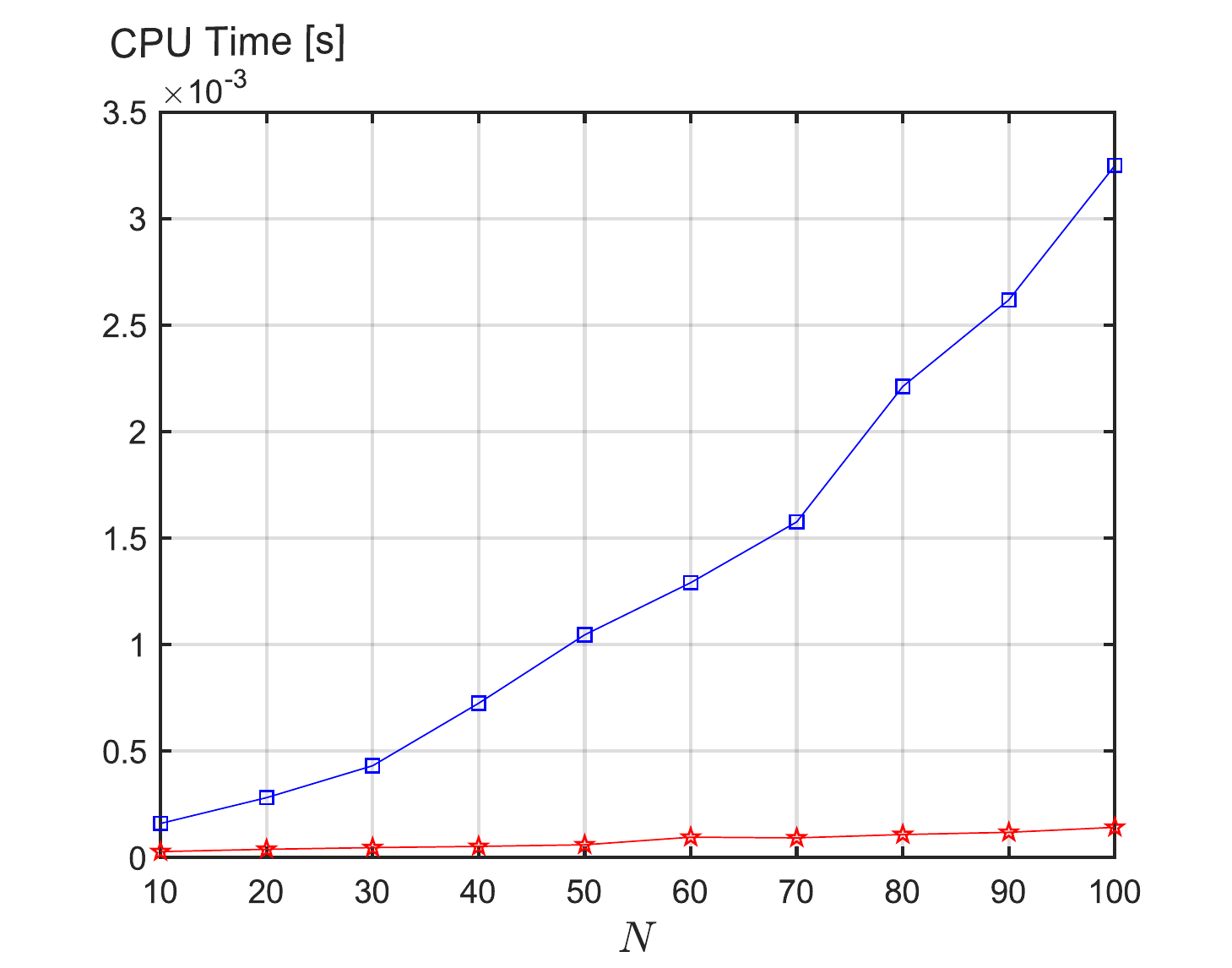}  
	\caption{Online computational time for the rigid Tube MPC controller as a 
		function of $N$ for {\tt qpOASES} (blue line with square markers) and Algorithm~1 with $\overline{m}=5$ (red line 
		with star markers).}  
	\label{fig::timecost} 
\end{figure}
Finally, we compared the performance of Algorithm~1 with {\tt qpOASES}.
Figure~\ref{fig::timecost} shows this comparison in terms of CPU time vs
$N$. For short horizons, the CPU time is comparable (for $N=10$ {\tt qpOASES} requires $150 \, [\mu\text{s}]$, while Algorithm~1 requires $28 \, [\mu\text{s}]$). It is also clear that as $N$ increases,  the computational time for the online solver grows faster
than that of Algorithm~1. For example, for $N=100$ {\tt qpOASES} takes $3 \, [\text{ms}]$,
while Algorithm~1 needs $0.14\, [\text{ms}]$ per real-time iteration only.

%%%%%%%%%%%%%%%%%%%%%%%%%%%%%%%%%%%%%%%%%%%%%%%
\section{Conclusion}
\label{sec::conclusion}
This paper has introduced a real-time implementation of rigid Tube MPC for 
discrete-time linear systems with additive uncertainty and polytopic
state, control and uncertainty constraints. The implementation is based on a parallelizable MPC scheme and requires, at each time step, the evaluation of precomputed piecewise affine maps and a linearly constrained quadratic program. The 
Tube MPC problem is solved suboptimally, but the algorithm maintains guarantees in terms of recursive feasibility, robust constraint satisfaction, and robust asymptotic stability.
% The proposed implementation has the potential of fulfilling the computational time and memory
% requirements of embedded hardware.
The approach has been illustrated on 
a case study where comparisons were made with both Explicit MPC implementations
as well as auto generated online QP solvers.
% More large-scale benchmark case studies will be part of future work.

%%%%%%%%%%%%%%%%%%%%%%%%%%%%%%%%%%%%%%%%%%%%%%%%%%%%%
\section*{ACKNOWLEDGMENTS}

{\footnotesize
%K.W., Y.J., M.E.V, and B.H. acknowledge financial support via ShanghaiTech University, Grant-Nr. F-0203-14-012. 
%J.O.~thanks the Scientific Grant Agency of the Slovak Republic for finacial support through the grants 1/0112/16, 1/0585/19, and the Slovak Research and Development Agency under the project APVV-15-0007.
All authors were supported via a bilateral grant
between China and Slovakia (SK-CN-2017-PROJECT-6558, APVV
SK-CN-2017-0026), KW, YJ, MEV, and BH acknowledge financial support via ShanghaiTech
U.~Grant F-0203-14-012. JO thanks the Scientific
Grant Agency of the Slovak Republic for finacial support, grants 1/0112/16, 1/0585/19, and the Slovak Research
and Dev.~Agency, APVV-15-0007.
}

\addtolength{\textheight}{-3cm}   
% reference
\bibliographystyle{abbrv}	
\bibliography{main.bib}

\end{document}